\documentclass[10pt]{article}

\usepackage{amsmath,amsfonts,amssymb,amsthm,epsfig}

\def\C{\mathbb{C}}
\def\Z{\mathbb{Z}}

\def\id{\mathbf{1}}

\def\g{\ensuremath{\mathfrak{g}}}

\def\V{\mathbf{V}}
\def\W{\mathbf{W}}

\DeclareMathOperator{\Hom}{Hom}

\DeclareMathOperator{\End}{End}
\DeclareMathOperator{\Aut}{Aut}
\DeclareMathOperator{\inc}{in}
\DeclareMathOperator{\out}{out}
\DeclareMathOperator{\tr}{tr}
\DeclareMathOperator{\Irr}{Irr}

\newtheorem{theo}{Theorem}[section]
\newtheorem{prop}[theo]{Proposition}
\newtheorem{lem}[theo]{Lemma}

\newtheorem{defin}[theo]{Definition}

\numberwithin{equation}{section}

\newcommand{\Span}{\operatorname{Span}}

\begin{document}
\title{A geometric realization of spin representations and Young diagrams
from quiver varieties}
\author{Alistair Savage}
\date{}
\maketitle
\noindent

\begin{abstract}
Applying the techniques of an earlier paper with Frenkel, we
develop a geometric realization of spin representations and
Clifford algebras. In doing so, we give an explicit
parametrization of the irreducible components of Nakajima
varieties of type $D$ in terms of Young diagrams. We explicity
compute the geometric action of the Lie algebra and are able to
extend the geometric action to the entire Clifford algebra used in
the classical construction of the spin representations.
\end{abstract}

\renewcommand{\thefootnote}{\fnsymbol{footnote}}

\section*{Introduction}
In \cite{FS03}, we related two apparently different bases in the
representations of affine Lie algebras of type $A$: one arising
from statistical mechanics, the other from gauge theory.  In
particular, using geometric methods associated to quiver
varieties, we were able to give an alternative and much simpler
geometric proof of a result of \cite{D89} on the construction of
bases of affine Lie algebra representations. At the same time, we
gave a simple parametrization of the irreducible components of
Nakajima quiver varieties associated to infinite and cyclic
quivers in terms of Young and Maya diagrams.  In the current
paper, we consider the spin representations of the Lie algebra $\g
= \mathfrak{so}_{2n}\C$ of type $D_n$.  Applying the techniques of
\cite{FS03}, we are able to give a very explicit parametrization
of the irreducible components of the associated Nakajima varieties
in terms of Young diagrams of strictly decreasing row length with
maximum length $n-1$.  We also explicitly compute the geometric
action and thus obtain a realization of the spin representations
in terms of these diagrams.  Furthermore, we are able to extend
the geometric action to the entire Clifford algebra used in the
classical construction of the spin representations.  This is the
first example of a geometric realization of additional structure
on irreducible representations extending the module structure
originally defined by Nakajima. Possible extensions of these
results to the affine case may suggest new statistical mechanics
results for type $D$ analogous to those of type $A$ mentioned
above.

The author would like to thank I. B. Frenkel for numerous
discussions and helpful suggestions and M. Kassabov for pointing
out an error in an earlier version of this paper.  This research
was supported in part by the Natural Sciences and Engineering
Research Council (NSERC) of Canada and the Clay Mathematics
Institute.

\section{Lusztig's Quiver Varieties}
\label{sec:lus_def}

In this section, we will recount the explicit description given in
\cite{L91} of the
irreducible components of Lusztig's quiver variety for
type $D_n$.  See this reference for details, including proofs.

\subsection{General Definitions}

Let $I$ be the set of vertices of the Dynkin graph of a symmetric
Kac-Moody Lie algebra $\mathfrak{g}$ and let $H$ be the set of
pairs consisting of an edge together with an orientation of it.
For $h \in H$, let $\inc(h)$ (resp. $\out(h)$) be the incoming
(resp. outgoing) vertex of $h$.  We define the involution $\bar{\
}: H \to H$ to be the function which takes $h \in H$ to the
element of $H$ consisting of the same edge with opposite
orientation.  An \emph{orientation} of our graph is a choice of a
subset $\Omega \subset H$ such that $\Omega \cup \bar{\Omega} = H$
and $\Omega \cap \bar{\Omega} = \emptyset$.

Let $\mathcal{V}$ be the category of
finite-dimensional $I$-graded vector spaces $\V = \oplus_{i
  \in I} \V_i$ over $\C$ with morphisms being linear maps
respecting the grading.  Then $\V \in \mathcal{V}$ shall denote
that $\V$ is an object of $\mathcal{V}$.  Labelling the vertices of
$I$ by $1,\dots,n$, the dimension of
$\V \in \mathcal{V}$ is given by $\mathbf{v} = \dim \V
= (\dim \V_1, \dots,
\dim \V_n)$.

Given $\V \in \mathcal{V}$, let
\[
\mathbf{E_V} = \bigoplus_{h \in H} \Hom (\V_{\out(h)},
\V_{\inc(h)}).
\]
For any subset $H'$ of $H$, let $\mathbf{E}_{\V, H'}$ be the
subspace of $\mathbf{E_V}$ consisting of all vectors $x = (x_h)$ such
that $x_h=0$ whenever $h \not\in H'$.  The algebraic group $G_\V
= \prod_i \Aut(\V_i)$ acts on $\mathbf{E_V}$ and
$\mathbf{E}_{\V, H'}$ by
\[
(g,x) = ((g_i), (x_h)) \mapsto gx = (x'_h) = (g_{\inc(h)} x_h g_{\out(h)}^{-1}).
\]

Define the function $\varepsilon : H \to \{-1,1\}$ by $\varepsilon
(h) = 1$ for all $h \in \Omega$ and $\varepsilon(h) = -1$ for all
$h \in {\bar{\Omega}}$.  For $\V \in \mathcal{V}$, the Lie algebra
of $G_\V$ is $\mathbf{gl_V} = \prod_i \End(\V_i)$ and it acts on
$\mathbf{E_V}$ by
\[
(a,x) = ((a_i), (x_h)) \mapsto [a,x] = (x'_h) = (a_{\inc(h)}x_h - x_h
a_{\out(h)}).
\]
Let $\left<\cdot,\cdot\right>$ be the nondegenerate,
$G_\V$-invariant, symplectic form on
$\mathbf{E_V}$ with values in $\C$ defined by
\[
\left<x,y\right> = \sum_{h \in H} \varepsilon(h) \tr (x_h y_{\bar{h}}).
\]
Note that $\mathbf{E_V}$ can be considered as the cotangent space of
$\mathbf{E}_{\V, \Omega}$ under this form.

The moment map associated to the $G_{\mathbf{V}}$-action on the
symplectic vector space $\mathbf{E_V}$ is the map $\psi : \mathbf{E_V}
\to \mathbf{gl_V}$ with $i$-component $\psi_i : \mathbf{E_V} \to \End
\V_i$ given by
\[
\psi_i(x) = \sum_{h \in H,\, \inc(h)=i} \varepsilon(h) x_h x_{\bar{h}} .
\]

\begin{defin}[\cite{L91}]
\label{def:nilpotent}
An element $x \in \mathbf{E_V}$ is said to be \emph{nilpotent} if
there exists an $N \ge 1$ such that for any sequence $h'_1, h'_2,
\dots, h'_N$ in $H$ satisfying $\out (h'_1) = \inc (h'_2)$, $\out (h'_2) =
\inc (h'_3)$, \dots, $\out (h'_{N-1}) = \inc (h'_N)$, the composition
$x_{h'_1} x_{h'_2} \dots x_{h'_N} : \V_{\out (h'_N)} \to
  \V_{\inc (h'_1)}$ is zero.
\end{defin}

\begin{defin}[\cite{L91}] Let $\Lambda_\V$ be the set of all
  nilpotent elements $x \in \mathbf{E_V}$ such that $\psi_i(x) = 0$
  for all $i \in I$.
\end{defin}


\subsection{Type $D_n$}
\begin{figure}
\begin{center}
\epsfig{file=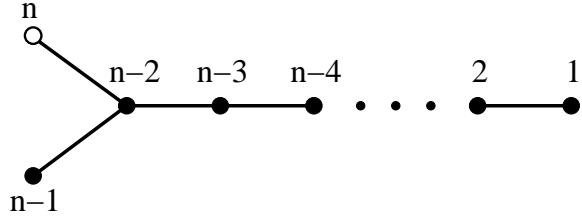,width=3in}
\caption{The Dynkin graph of type $D_n$. We represent the $n$th vertex by
an open dot to distinguish it from the $(n-1)$st vertex.
\label{fig:dynkin_dn}}
\end{center}
\end{figure}

\begin{figure}
\begin{center}
\epsfig{file=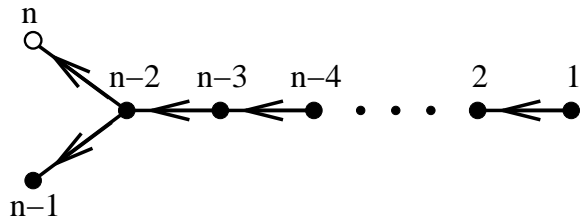,width=3in}
\caption{The quiver (oriented graph) of type $D_n$. \label{fig:quiver_dn}}
\end{center}
\end{figure}

Let $\g=\mathfrak{so}_{2n}\C$ be the simple Lie algebra of type
$D_n$. Let $I=\{1,2,\dots,n\}$ be the set of vertices of the
Dynkin graph of \g, labelled as in Figure~\ref{fig:dynkin_dn}.
We let $\Omega$ be the orientation indicated in Figure~\ref{fig:quiver_dn}.

Label each oriented edge by its incoming and outgoing vertices.
That is, if vertices $i$ and $j$ are connected by an edge, $h_{i,j}$
denotes the oriented edge with $\out(h)=i$, $\inc(h)=j$.
The following is proven in \cite[Proposition 14.2]{L91}.

\begin{prop}
\label{prop:irrcomp:Dn}
The irreducible components of $\Lambda_\V$ are the closures of
the conormal bundles of the various $G_\V$-orbits in
$\mathbf{E}_{\V, \Omega}$.
\end{prop}

\begin{figure}
\begin{center}
\epsfig{file=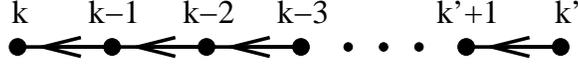,width=3in}
\caption{The string representing $(\V(k',k),x(k',k))$, $1 \le k' \le k \le
n-1$. \label{fig:lowerstring}}
\end{center}
\end{figure}

For two integers $1 \le k' \le k \le n-1$, define $\V(k',k)
\in \mathcal{V}$
to be the vector space with basis $\{ e_r\ |\ k' \le r \le k\}$.  We
require that $e_r$ has degree $r \in I$.
Let $x(k', k) \in \mathbf{E}_{\V(k',k), \Omega}$
be defined by
$x(k',k)_{h_{r,r+1}} : e_r \mapsto e_{r+1}$ for $k' \le r < k$, and
all other components of $x(k',k)$ are zero.
We picture this representation as the string of Figure~\ref{fig:lowerstring}.

\begin{figure}
\begin{center}
\epsfig{file=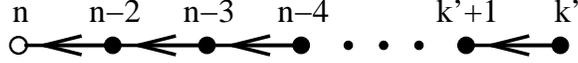,width=3in}
\caption{The string representing $(\V(k',n),x(k',n))$, $1 \le k' \le n-2$
or $k'=n$. \label{fig:upperstring}}
\end{center}
\end{figure}

For an integer $k'$ such that $1 \le k' \le n-2$ or $k'=n$, define $\V(k',n)$ to
be the vector space with basis $\{e_r\ |\ k' \le r \le n-2 \text{ or }
r=n\}$.  Again, we require that $e_r$ has degree $r \in I$.  Let
$x(k',n) \in \mathbf{E}_{\V(k',n), \Omega}$ be defined by
\begin{align*}
x(k',n)_{h_{r,r+1}}(e_r) &= \begin{cases}
  e_{r+1} & \text{if $r < n-2$} \\
  0 & \text{otherwise}
\end{cases},\\
x(k',n)_{h_{n-2,n}}(e_{n-2}) &= e_n
\end{align*}
and all other components of $x(k',k)$ are zero.
We picture this representation as the string of Figure~\ref{fig:upperstring}.

\begin{figure}
\begin{center}
\epsfig{file=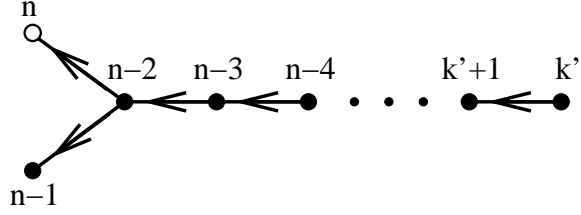,width=3in}
\caption{The (forked) string representing $(\V(k',n+1),x(k',n+1))$,
$1 \le k' \le n-2$. \label{fig:forkstring}}
\end{center}
\end{figure}

Next, for an integer $k'$ such that $1 \le k' \le n-2$, define
$\V(k',n+1)$ to be the vector space with basis $\{e_r\ |\ k' \le r
\le n\}$ where the degree of $e_r$ is $r \in I$.  $x(k',n+1)$ is
defined by
\begin{align*}
x(k',n+1)_{h_{r,r+1}}(e_r) &= \begin{cases}
  e_{r+1} & \text{if $r \le n-2$} \\
  0 & \text{otherwise}
\end{cases},\\
x(k',n+1)_{h_{n-2,n}}(e_{n-2}) &= e_n
\end{align*}
and all other components of $x(k',n+1)$ are zero.
We picture this representation as the (forked) string of
Figure~\ref{fig:forkstring}.

\begin{figure}
\begin{center}
\epsfig{file=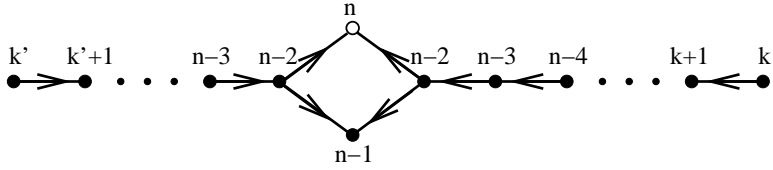,width=4in} \caption{The string
representing $({\widetilde \V}(k',k),{\tilde x}(k',k))$, $1 \le k'
< k \le n-2$. \label{fig:doublestring}}
\end{center}
\end{figure}

Finally, for $1 \le k' < k \le n-2$, let ${\widetilde \V}(k',k)$
be the vector space with basis $\{e_r\, |\, k \le r \le n\} \cup
\{{\tilde e}_r\, |\, k' \le r \le n-2\}$ where the degree of $e_r$
and ${\tilde e}_r$ is $r \in I$.  Let ${\tilde x}(k',k)$ be
defined by
\begin{align*}
{\tilde x}(k',k)_{h_{r,r+1}}(e_r) &= \begin{cases} e_{r+1} &
\text{if $r \ge n-2$} , \\
0 & \text{otherwise} \end{cases}, \\
{\tilde x}(k',k)_{h_{r,r+1}}({\tilde e}_r) &= \begin{cases}
{\tilde e}_{r+1} & \text{if $r < n-2$} , \\
e_{n-1} & \text{if $r=n-2$} \\
0 & \text{otherwise} \end{cases}, \\
{\tilde x}(k',k)_{h_{n-2,n}}(e_{n-2}) &= e_n \\
{\tilde x}(k',k)_{h_{n-2,n}}({\tilde e}_{n-2}) &= e_n
\end{align*}
and all other components of ${\tilde x}(k',k)$ are zero.  We
picture this representation as the string of
Figure~\ref{fig:doublestring}.

\begin{prop}
\label{prop:ind-reps} The above $(\V(k',k), x(k',k))$ and
$({\widetilde \V}(k',k), {\tilde x}(k',k))$ are indecomposable
representations of the $D_n$ quiver with orientation $\Omega$.
Conversely, any indecomposable finite-dimensional representation
$(\V,x)$ of this quiver is isomorphic to one of these
representations.
\end{prop}

\begin{proof}
It is easy to see that each of the above representations is
indecomposable.  The fact that our list is exhaustive follows from
the fact that the indecomposable representations of the $D_n$
quiver are in one-to-one correspondence with the positive roots of
the Lie algebra of type $D_n$ (see \cite{K80,K83}).
\end{proof}

Let
\begin{align*}
Z =& \{(k',k)\, |\, (k',k) \ne (n-1,n),(n-1,n+1),(n,n+1),
(n+1,n+1)\} \\
&\cup \{(k',k)^{\sim} \, |\, 1 \le k' < k \le n-2\},
\end{align*}
and let $\tilde Z$ be the set of all functions $Z \to \Z_{\ge 0}$
with finite support. It is clear that for $\V \in \mathcal{V}$,
the set of $G_\V$-orbits in $\mathbf{E}_{\V, \Omega}$ is naturally
indexed by the subset $\tilde Z_\V$ of $\tilde Z$ consisting of
those $f \in \tilde Z$ such that
\[
\sum_{A(k',k) \ni i} f(k',k) + \sum_{B(k',k) \ni i} f((k',k)^\sim)
(\dim {\widetilde \V}(k',k))_i = \dim \V_i \quad \forall \ i \in I
\]
where $A(k',k)$ is the set of all $i \in I$ such that $\V(k',k)$
contains a vector of degree $i$ and $B(k',k)$ is the set of all $i
\in I$ such that ${\widetilde \V}(k',k)$ contains a vector of
degree $i$. Note that we write $f(k',k)$ for $f((k',k))$.
Corresponding to a given $f$ is the orbit consisting of all
representations isomorphic to a sum of the indecomposable
representations $x(k',k)$ and ${\tilde x}(k',k)$, each occurring
with multiplicity $f(k',k)$ and $f((k',k)^\sim)$ respectively.
Denote by $\mathcal{O}_f$ the $G_\V$-orbit corresponding to $f \in
\tilde Z_\V$.  Let $\mathcal{C}_f$ be the conormal bundle to
$\mathcal{O}_f$ and let $\bar{\mathcal{C}}_f$ be its closure.  We
then have the following proposition.

\begin{prop}
The map $f \to \bar{\mathcal{C}}_f$ is a one-to-one correspondence
between the set
$\tilde Z_\V$ and the set of irreducible
components of $\Lambda_\V$.
\end{prop}
\begin{proof}
This follows immediately from Propositions~\ref{prop:irrcomp:Dn}
and \ref{prop:ind-reps}.
\end{proof}


\section{Nakajima's Quiver Varieties}
\label{sec:def_nak}

We introduce here a description of the quiver varieties first
presented in \cite{N94} for type $D_n$.

\begin{defin}[\cite{N94}]
\label{def:lambda}
For $\mathbf{v}, \mathbf{w} \in \Z_{\ge 0}^I$, choose $I$-graded
vector spaces $\V$ and $\W$ of graded dimensions
$\mathbf{v}$ and
$\mathbf{w}$ respectively.  Then define
\[
\Lambda \equiv \Lambda(\mathbf{v},\mathbf{w}) =
\Lambda_\V \times \sum_{i \in I} \Hom (\V_i, \W_i).
\]
\end{defin}

Now, suppose that $\mathbf{S}$ is an $I$-graded subspace of $\V$.
For $x \in
\Lambda_\mathbf{V}$ we say that $\mathbf{S}$ is
\emph{$x$-stable} if $x(\mathbf{S}) \subset \mathbf{S}$.

\begin{defin}[\cite{N94}]
\label{def:lambda-stable}
Let $\Lambda^{\text{st}} = \Lambda(\mathbf{v},\mathbf{w})^{\text{st}}$ be
the set of all $(x, j) \in
\Lambda(\mathbf{v},\mathbf{w})$ satisfying the following condition:  If
$\mathbf{S}=(\mathbf{S}_i)$ with $\mathbf{S}_i \subset \mathbf{V}_i$ is
$x$-stable and $j_i(\mathbf{S}_i) = 0$ for
$i \in I$, then $\mathbf{S}_i = 0$ for $i \in I$.
\end{defin}

The group $G_\mathbf{V}$ acts on $\Lambda(\mathbf{v},\mathbf{w})$ via
\[
(g,(x,j)) = ((g_i), ((x_h), (j_i))) \mapsto ((g_{\inc (h)} x_h
g_{\out (h)}^{-1}), (j_i g_i^{-1})).
\]
and the stabilizer of any point of
$\Lambda(\mathbf{v},\mathbf{w})^{\text{st}}$ in $G_{\mathbf{V}}$ is trivial
  (see \cite[Lemma~3.10]{N98}).  We then make the following definition.
\begin{defin}[\cite{N94}]
\label{def:L}
Let $\mathcal{L} \equiv \mathcal{L}(\mathbf{v},\mathbf{w}) =
\Lambda(\mathbf{v},\mathbf{w})^{\text{st}} / G_{\mathbf{V}}$.
\end{defin}

Let $\Irr \mathcal{L}(\mathbf{v},\mathbf{w})$ (resp. $\Irr
\Lambda(\mathbf{v},\mathbf{w})$) be the set of irreducible components of
$\mathcal{L}(\mathbf{v},\mathbf{w})$ (resp. $\Lambda(\mathbf{v},\mathbf{w})$).
Then $\Irr \mathcal{L}(\mathbf{v},\mathbf{w})$ can be identified with
\[
\{ Y \in \Irr \Lambda(\mathbf{v},\mathbf{w})\, |\, Y \cap
\Lambda(\mathbf{v},\mathbf{w})^{\text{st}} \ne \emptyset \}.
\]
Specifically, the irreducible components of $\Irr
\mathcal{L}(\mathbf{v},\mathbf{w})$ are precisely those
\[
X_f \stackrel{\text{def}}{=} \left( \left( \bar{\mathcal{C}}_f \times
    \sum_{i \in I} \Hom (\mathbf{V}_i, \mathbf{W}_i) \right) \cap
\Lambda(\mathbf{v},\mathbf{w})^{\text{st}} \right) / G_\mathbf{V}
\]
which are nonempty.

The following will be used in the sequel.
\begin{lem}
\label{lem:irrcomp}
One has
\[
\Lambda^{\text{st}} = \left\{ x \in \Lambda\, \left| \,
       \ker j_i \cap \bigcap_{h: \out(h)=i} \ker x_h  = 0
      \ \forall \ i\right. \right\} .
\]
\end{lem}
\begin{proof}
Since each $\bigcap_{h: \out{h}=i} \ker x_h$ is $x$-stable, the
left hand side is obviously contained in the right hand side.  Now
suppose $x$ is an element of the right hand side.  Let $\mathbf{S}
= (\mathbf{S}_i)$ with $\mathbf{S}_i \subset \mathbf{V}_i$ be
$x$-stable and $j_i (\mathbf{S}_i) = 0$ for $i \in I$.  Assume
that $\mathbf{S} \ne 0$. Since all elements of $\Lambda$ are
nilpotent, we can find a minimal value of $N$ such that the
condition in Definition~\ref{def:nilpotent} is satisfied.  Then we
can find a $v \in \mathbf{S}_i$ for some $i$ and a sequence $h_1',
h_2', \dots, h_{N-1}'$ (empty if $N=1$) in $H$ such that $\out
(h_1') = \inc (h_2')$, $\out (h_2') = \inc (h_3')$, \dots, $\out
(h_{N-2}') = \inc (h_{N-1}')$ and $v' = x_{h_1'} x_{h_2'} \dots
x_{h_{N-1}'} (v) \ne 0$. Now, $v' \in \mathbf{S}_{i'}$ for some
$i' \in I$ by the stability of $\mathbf{S}$ (hence $j_{i'}(v') =
0$) and $v' \in \bigcap_{h: \out{h}=i'} \ker x_h$ by our choice of
$N$.  This contradicts the fact that $x$ is an element of the
right hand side.
\end{proof}


\section{The Lie Algebra Action}
\label{sec:structure}

We summarize here some results from \cite{N94} that will be needed in
the sequel.  See this reference for more details, including proofs.
We keep the notation of Sections~\ref{sec:lus_def}
and~\ref{sec:def_nak}.

Let $\mathbf{w, v, v', v''} \in \Z_{\ge 0}^I$ be such that $\mathbf{v}
= \mathbf{v'} + \mathbf{v''}$.  Consider the maps
\begin{equation}
\label{eq:diag_action}
\Lambda(\mathbf{v}'',\mathbf{0}) \times \Lambda(\mathbf{v}',\mathbf{w})
\stackrel{p_1}{\leftarrow} \mathbf{\tilde F (v,w;v'')}
\stackrel{p_2}{\rightarrow} \mathbf{F(v,w;v'')}
\stackrel{p_3}{\rightarrow} \Lambda(\mathbf{v},\mathbf{w}),
\end{equation}
where the notation is as follows.  A point of $\mathbf{F(v,w;v'')}$ is
a point $(x,j) \in \Lambda(\mathbf{v},\mathbf{w})$ together with an $I$-graded,
$x$-stable
subspace $\mathbf{S}$ of $\mathbf{V}$ such that $\dim \mathbf{S} =
\mathbf{v'} = \mathbf{v} - \mathbf{v''}$.  A point of $\mathbf{\tilde
  F (v,w;v'')}$ is a point $(x,j,\mathbf{S})$ of $\mathbf{F(v,w;v'')}$
together with a collection of isomorphisms $R'_i : \mathbf{V}'_i \cong
\mathbf{S}_i$ and $R''_i : \mathbf{V}''_i \cong \mathbf{V}_i /
\mathbf{S}_i$ for each $i \in I$.  Then we define $p_2(x,j,\mathbf{S},
R',R'') = (x,j,\mathbf{S})$, $p_3(x,j,\mathbf{S}) = (x,j)$ and
$p_1(x,j,\mathbf{S},R',R'') = (x'',x',j')$ where $x'', x', j'$ are
determined by
\begin{align*}
R'_{\inc(h)} x'_h &= x_h R'_{\out(h)} : \mathbf{V}'_{\out(h)} \to
\mathbf{S}_{\inc(h)}, \\
j'_i &= j_i R'_i : \mathbf{V}'_i \to \mathbf{W}_i \\
R''_{\inc(h)} x''_h &= x_h R''_{\out(h)} : \mathbf{V}''_{\out(h)} \to
\mathbf{V}_{\inc(h)} / \mathbf{S}_{\inc(h)}.
\end{align*}
It follows that $x'$ and $x''$ are nilpotent.

\begin{lem}[{\cite[Lemma 10.3]{N94}}]
One has
\[
(p_3 \circ p_2)^{-1} (\Lambda(\mathbf{v},\mathbf{w})^{\text{st}}) \subset
p_1^{-1} (\Lambda(\mathbf{v}'',\mathbf{0}) \times
\Lambda(\mathbf{v}',\mathbf{w})^{\text{st}}).
\]
\end{lem}

Thus, we can restrict \eqref{eq:diag_action} to
$\Lambda^{\text{st}}$, forget the
$\Lambda(\mathbf{v}'',\mathbf{0})$-factor and consider the quotient by
$G_\mathbf{V}$, $G_\mathbf{V'}$.  This yields the diagram
\begin{equation}
\label{eq:diag_action_mod}
\mathcal{L}(\mathbf{v'}, \mathbf{w}) \stackrel{\pi_1}{\leftarrow}
\mathcal{F}(\mathbf{v}, \mathbf{w}; \mathbf{v} - \mathbf{v'})
\stackrel{\pi_2}{\rightarrow} \mathcal{L}(\mathbf{v}, \mathbf{w}),
\end{equation}
where
\[
\mathcal{F}(\mathbf{v}, \mathbf{w}, \mathbf{v} - \mathbf{v'})
\stackrel{\text{def}}{=} \{ (x,j,\mathbf{S}) \in \mathbf{F(v,w;v-v')}\,
  |\, (x,j) \in \Lambda(\mathbf{v},\mathbf{w})^{\text{st}} \} / G_\mathbf{V}.
\]

Let $M(\mathcal{L}(\mathbf{v}, \mathbf{w}))$ be the vector space of
all constructible functions on $\mathcal{L}(\mathbf{v}, \mathbf{w})$.
For a subvariety $Y$ of a variety $A$, let $\mathbf{1}_Y$ denote the
function on $A$ which takes the value 1 on $Y$ and 0 elsewhere.  Let
$\chi (Y)$ denote the Euler characteristic of the algebraic variety
$Y$.  Then for a map $\pi$ between algebraic varieties $A$ and $B$, let
$\pi_!$ denote the map between the abelian groups of constructible
functions on $A$ and $B$ given by
\[
\pi_! (\mathbf{1}_Y)(y) = \chi (\pi^{-1}(y) \cap Y),\ Y \subset A
\]
and let $\pi^*$ be the pullback map from functions on $B$ to functions
on $A$ acting as $\pi^* f(y) = f(\pi(y))$.
Then define
\begin{align*}
&H_i : M(\mathcal{L}(\mathbf{v}, \mathbf{w})) \to
M(\mathcal{L}(\mathbf{v}, \mathbf{w})); \quad H_i f = u_i f, \\
&E_i : M(\mathcal{L}(\mathbf{v}, \mathbf{w})) \to
M(\mathcal{L}(\mathbf{v} - \mathbf{e}^i, \mathbf{w})); \quad E_i f =
(\pi_1)_! (\pi_2^* f), \\
&F_i : M(\mathcal{L}(\mathbf{v} - \mathbf{e}^i, \mathbf{w})) \to
M(\mathcal{L}(\mathbf{v}, \mathbf{w})); \quad F_i g = (\pi_2)_!
(\pi_1^* g).
\end{align*}
Here
\[
\mathbf{u} = {^t(u_0, \dots, u_n)} = \mathbf{w} - C \mathbf{v}
\]
where $C$ is the Cartan matrix of $\mathfrak{g}$ and we are using
diagram~\eqref{eq:diag_action_mod} with $\mathbf{v}' = \mathbf{v}
- \mathbf{e}^i$ where $\mathbf{e}^i$ is the vector whose
components are given by $\mathbf{e}^i_j = \delta_{ij}$.

Now let $\varphi$ be the constant function on $\mathcal{L}(\mathbf{0},
\mathbf{w})$ with value 1.  Let $L(\mathbf{w})$ be the vector space of
functions generated by acting on $\varphi$ with all possible combinations of
the operators $F_i$.  Then let $L(\mathbf{v},\mathbf{w}) =
M(\mathcal{L}(\mathbf{v}, \mathbf{w})) \cap L(\mathbf{w})$.

\begin{prop}[{\cite[Thm 10.14]{N94}}]
The operators $E_i$, $F_i$, $H_i$ on $L(\mathbf{w})$ provide the
structure of the
irreducible highest weight integrable representation of
$\mathfrak{g}$ with highest weight $\mathbf{w}$.  Each
summand of the decomposition $L(\mathbf{w}) = \bigoplus_\mathbf{v}
L(\mathbf{v}, \mathbf{w})$ is a weight space with weight
$\mathbf{w} - C \mathbf{v}$.  That is, with weight
\[
\sum_{i \in I} (\mathbf{w} - C\mathbf{v})_i \Lambda_i
\]
where the $\Lambda_i$ are the fundamental weights of \g\ (i.e.
$\Lambda_i (\alpha_j) = \delta_{ij}$ where $\alpha_j$ is the
simple root corresponding to the vertex $j$).
\end{prop}

Let $X \in \Irr \mathcal{L}(\mathbf{v}, \mathbf{w})$ and define a
linear map $T_X : L(\mathbf{v}, \mathbf{w}) \to \C$ as in
\cite[3.8]{L92}.  The map $T_X$ associates to a constructible function $f \in
L(\mathbf{v}, \mathbf{w})$ the (constant) value of $f$ on a suitable
open dense subset of $X$.  The fact that $L(\mathbf{v}, \mathbf{w})$
is finite-dimensional allows us to take such an open set on which
\emph{any} $f \in L(\mathbf{v}, \mathbf{w})$ is constant.  So we have
a linear map
\[
\Phi : L(\mathbf{v}, \mathbf{w}) \to \C^{\Irr \mathcal{L}(\mathbf{v},
  \mathbf{w})}.
\]
The following proposition is proved in \cite[4.16]{L92} (slightly
generalized in
\cite[Proposition 10.15]{N94}).

\begin{prop}
\label{prop:func_irrcomp_isom}
The map $\Phi$ is an isomorphism; for any $X \in \Irr \mathcal{L}(\mathbf{v},
\mathbf{w})$, there is a unique function $g_X \in L(\mathbf{v},
\mathbf{w})$ such that for some open dense subset $O$ of $X$ we have
$g_X|_O = 1$ and for some closed $G_\mathbf{V}$-invariant
subset $K \subset \mathcal{L}(\mathbf{v}, \mathbf{w})$ of dimension $<
\dim \mathcal{L}(\mathbf{v}, \mathbf{w})$ we have $g_X=0$ outside $X
\cup K$.  The functions $g_X$ for $X \in \Irr \Lambda(\mathbf{v},\mathbf{w})$
form a basis of $L(\mathbf{v},\mathbf{w})$.
\end{prop}


\section{Geometric Realization of the Spin Representations}
\label{sec:level1}

We now seek to describe the irreducible components of Nakajima's
quiver variety corresponding to the spin representations of the Lie
algebra \g\ of type $D_n$.  By the
comment made in
Section~\ref{sec:def_nak}, it
suffices to determine which irreducible components of
$\Lambda(\mathbf{v},\mathbf{w})$ are not killed by the stability
condition.  By
Definition~\ref{def:lambda} and Lemma~\ref{lem:irrcomp}, these are
precisely those irreducible components which contain points $x$ such
that
\begin{equation}
\label{eq:killcond}
\dim \left( \bigcap_{h: \out(h)=i} \ker x_h \right)
  \le \mathbf{w}_i \ \forall\ i.
\end{equation}

We will consider the spin representations of highest weights
$\Lambda_{n-1}$ and $\Lambda_n$ which correspond to $\mathbf{w} =
\mathbf{w}^{n-1}$ and $\mathbf{w}^n$ respectively, where
$\mathbf{w}^i_j = \delta_{ij}$.

Let $\mathcal{Y}$ be the set of all Young diagrams with strictly
decreasing row lengths of length at most $n-1$, that is, the set
of all strictly decreasing sequences $(l_1, \dots, l_s)$ of
non-negative integers ($l_j=0$ for $j>s$) such that $l_1 \le n-1$.
We will use the terms Young diagram and partition interchangeably.
For $Y=(l_1, \dots, l_s) \in \mathcal{Y}$ and $1 \le i \le s$, let
\begin{align*}
A_{l_i}^+ &= \begin{cases}
  (n-l_i,n) & \text{if $i$ is odd and $l_i > 1$} \\
  (n,n) & \text{if $i$ is odd and $l_i=1$} \\
  (n-l_i,n-1) & \text{if $i$ is even}
\end{cases}, \\
A_{l_i}^- &= \begin{cases}
  (n-l_i,n) & \text{if $i$ is even and $l_i > 1$} \\
  (n,n) & \text{if $i$ is even and $l_i=1$} \\
  (n-l_i,n-1) & \text{if $i$ is odd}
\end{cases}
\end{align*}
and let
\[
A_Y^{\pm} = \bigcup_{i=1}^s A_{l_i}^{\pm}.
\]

\begin{theo}
\label{thm:irrcomp_lev1_infty}
The irreducible components of $\mathcal{L}(\mathbf{v},\mathbf{w}^n)$
(resp. $\mathcal{L}(\mathbf{v},\mathbf{w}^{n-1})$)
are precisely those $X_f$ where $f \in \tilde Z_\V$ such that
\[
\{(k',k)\ |\ f(k',k)=1\} = A_Y^\pm
\]
for some $Y \in \mathcal{Y}$, $f(k',k)=0$ for $(k',k) \not \in
A_Y^\pm$, and $f((k',k)^\sim) = 0$ for all $(k',k)$. Denote the
component corresponding to such an $f$ by $X_Y^\pm$.  Thus, $Y
\leftrightarrow X_Y^\pm$ is a natural 1-1 correspondence between
the set $\mathcal{Y}$ and the irreducible components of
$\cup_\mathbf{v} \mathcal{L}(\mathbf{v}, \mathbf{w}^n)$ (resp.
$\cup_\mathbf{v} \mathcal{L}(\mathbf{v}, \mathbf{w}^{n-1})$). Here
the plus signs correspond to the case $\cup_\mathbf{v}
\mathcal{L}(\mathbf{v},\mathbf{w}^n)$ and the minus signs
correspond to the case $\cup_\mathbf{v} \mathcal{L}(\mathbf{v},
\mathbf{w}^{n-1})$.
\end{theo}

\begin{proof}
We prove only the $\cup_\mathbf{v}
\mathcal{L}(\mathbf{v},\mathbf{w}^n)$ case. The other case in
analogous. Consider the two representations
$(\V(k_1',k_1),x(k_1',k_1))$ and $(\V(k_2',k_2),x(k_2',k_2))$, $1
\le k' \le k \le n-1$ of our oriented graph as described in
Section~\ref{sec:lus_def} where the basis of $\V(k_i',k_i)$ is
$\{e^i_r\ |\ k_i' \le r \le k_i\}$. Let $W$ be the conormal bundle
to the $G_\V$-orbit through the point
\[
x_\Omega = (x_h)_{h \in \Omega} = x(k_1',k_1) \oplus x(k_2',k_2) \in
  \mathbf{E}_{\V(k_1',k_1) \oplus \V(k_2',k_2), \Omega}.
\]
By the proof of Theorem~5.1 of \cite{FS03}, we see that for
\[
x = (x_\Omega, x_{\bar \Omega}) = ((x_h)_{h \in \Omega}, (x_h)_{h \in
  {\bar \Omega}}) = (x_h)_{h \in H} \in W,
\]
all of
$V(k_2',k_2)$ must be in the kernel of $x_{\bar \Omega}$ unless $k_2'
> k_1'$ and $k_2 > k_1$.  If these conditions are satisfied, there
exists a point $(x_\Omega, x_{\bar \Omega}) \in W$ such that
$x_{h_{i,i-1}} (e^2_i) = c e^1_{i-1}$ for all $k_1+1 \le i \le
k_2'$ (for a fixed non-zero $c$). We picture such a representation
as in Figure~\ref{fig:youngcomm}.
\begin{figure}
\begin{center}
\epsfig{file=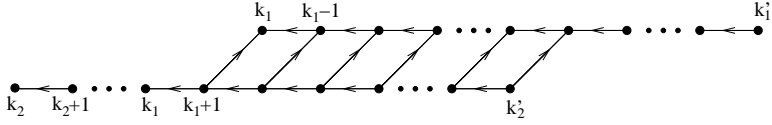,width=4in}
\caption{If $x_{h_{i,i-1}} (e^2_i) \ne 0$ for some $i$, the
  commutativity of
  the above diagram forces $k_2' > k_1'$ and
  $k_2 > k_1$. The vertices in the upper (resp. lower) string represent the
  basis vectors defining the
  representation $(\V(k_1',k_1),x(k_1',k_1))$
  (resp. $(\V(k_2',k_2),x(k_2',k_2))$).
  A vertex labelled $i$ represents $e^j_i \in \mathbf{V}_i$ ($j=1,2$).
  The arrows indicate the
  action of the obvious component of $x$.
 \label{fig:youngcomm}}
\end{center}
\end{figure}

Similarly, for the two represenations $(\V(k_1',n), x(k_1',n))$ and
$(\V(k_2',n), x(k_2',n))$, Both $\V(k_1',n)$ and $\V(k_2',n)$ are in the
kernel of any
component $x_h$, $h \in {\bar \Omega}$ of a point in the conormal
bundle to the orbit through $x(k_1',n) \oplus x(k_2',n)$ (since the left
endpoints of the strings are the same).

However, for $k_1' < k_2'$, there are points $x$ in the conormal bundle to the
orbit through $x(k_1',n) \oplus x(k_2',n-1)$ such that
$x_{h_{i,i-1}} (e^2_i) = c e^1_{i-1}$ for all $k_2' \le i \le n-1$ and some
non-zero $c$.  These maps do not violate the moment map condition since the
left endpoints of the two strings are different.  Such a representation is
pictured as in Figure~\ref{fig:keycomm}.  By symmetry, reversing the roles of
$n-1$ and $n$, there are points $x$ in the conormal bundle to the
orbit through $x(k_1',n-1) \oplus x(k_2',n)$ such that
$x_{h_{n,n-2}} (e^2_n) = c e^1_{n-2}$ and
$x_{h_{i,i-1}} (e^2_i) = c e^1_{i-1}$ for all $k_2' \le i \le n-2$ and some
non-zero $c$ if and only if $k_1' < k_2'$
\begin{figure}
\begin{center}
\epsfig{file=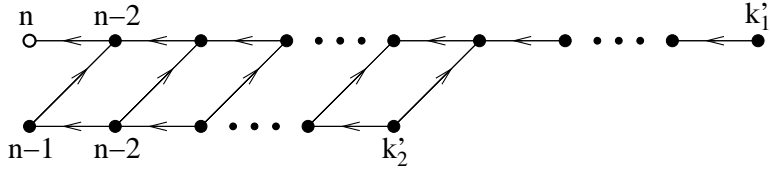,width=4in} \caption{If $x_{h_{i,i-1}}
(e^2_i) \ne 0$ for some $i$, the
  commutativity of
  the above diagram forces $k_2' > k_1'$. The vertices in the upper (resp. lower) string represent the
  basis vectors defining the
  representation $(\V(k_1',n),x(k_1',n))$
  (resp. $(\V(k_2',n-1),x(k_2',n-1))$).
  A vertex labelled $i$ represents $e^j_i \in \mathbf{V}_i$ ($j=1,2$).
  The arrows indicate the
  action of the obvious component of $x$.
 \label{fig:keycomm}}
\end{center}
\end{figure}

Now consider the representations $(\V(k_1',n-1),x(k_1',n-1))$ and
$(\V(k_2',n+1), x(k_2',n+1))$.  Suppose there exists a point
$(x_\Omega, x_{\bar \Omega})$ in the conormal bundle to the orbit
through the point $x(k_1',n-1) \oplus x(k_2',n+1)$ such that
$(\V(k_2',n+1), x(k_2',n+1))$ is not contained in the kernel of
$x_{\bar \Omega}$.  Then there is some basis element $e^2_i \in
\V(k_2',n+1)$ that is not killed by $x_{\bar \Omega}$.  This
implies that $x_{h_{i,i-1}}(e^2_i) = c e^1_{i-1}$ for some $c \ne
0$ since $x_{h_{i,i-1}}(e^2_i)$ can have no $e^2_{i-1}$ component
by nilpotency.  As the cases above, this implies $k_2' > k_1'$.

Now, let $x = (x_\Omega, x_{\bar \Omega})$ lie in the conormal bundle
to the point
\begin{equation}
\bigoplus_{i=1}^s x(k_i',k_i) \in \mathbf{E}_{\oplus_{i=1}^s
  \V(k_i', k_i), \Omega}.
\end{equation}
We can assume (by reordering the indices if necessary) that $k_1'
\le k_2' \le \dots \le k_s'$.  By the above arguments, if $k_1 =
n+1$ then both $e^1_{n-1}$ and $e^1_n$ lie in the kernel of
$x_{\bar \Omega}$ (and hence $x$).  But this violates the
stability condition.  So $k_1 \le n$.  Then we know from the above
that $e^1_{k_1}$ is in the kernel of $x_{\bar \Omega}$ (and hence
$x$) since $k_j' \ge k_1'$ for all $j$.  Thus
\[
e^1_{k_1} \in \bigcap_{h: \out(h)=k_1} \ker x_h .
\]

By the stability condition, we must then have $k_1=n$ and there can be
no other $e^j_i$ in $\bigcap_{h: \out(h)=i} \ker x_h$ for any $i$.
Suppose $k_2=n+1$.  Then we must have $x_{h_{n,n-2}}(e^2_n) =
ce^1_{n-2}$.  But then $x_{h_{n-2,n}} x_{h_{n,n-2}}(e^2_n) = ce^1_n
\ne 0$ which violates the moment map condition at the $n$th vertex.
So $k_2 \le n$.
Then, by the above considerations, $e^2_{k_2}$ is in
$\bigcap_{h: \out(h)=k_2} \ker x_h$ unless
$k_2=n-1$ and $x_{h_{n-1,n-2}}(e^2_{n-1})$ is a
non-zero multiple of $e^1_{n-2}$.  Continuing in this manner, we see
that we must have
\[
k_i = \begin{cases}
  n & \text{if $i$ is odd} \\
  n-1 & \text{if $i$ is even}
\end{cases}
\]
and $x_{h_{k_{i+1},n-2}}
(e^{i+1}_{k_{i+1}})$ has non-zero $e^i_{n-2}$ component for $1 \le i \le s-1$.
Then by the above we must have $k_{i+1}' > k_i'$ for $1 \le i \le s-1$.
Setting
\[
l_i = \begin{cases}
 n - k_i' & \text{if $k_i' \le n-2$} \\
 1       & \text{if $k_i'=n-1,n$}
\end{cases}
\]
we have an irreducible component of the type mentioned in the
theorem.  It remains to consider the case where a summand ${\tilde
x}(k',k)$ occurs with some non-zero multiplicity.  Using the
moment map condition and the same sort of arguments as above, one
can see that $x_h$ for $h \in {\bar \Omega}$ acting on vectors of
such summands can only have non-zero components of vectors in
other summands of this type or summands $x(k',k)$ for $k \le n-2$.
Thus, it follows by the nilpotency and stability conditions that
no such summands can occur.  The theorem follows.
\end{proof}

\begin{figure}
\begin{center}
\epsfig{file=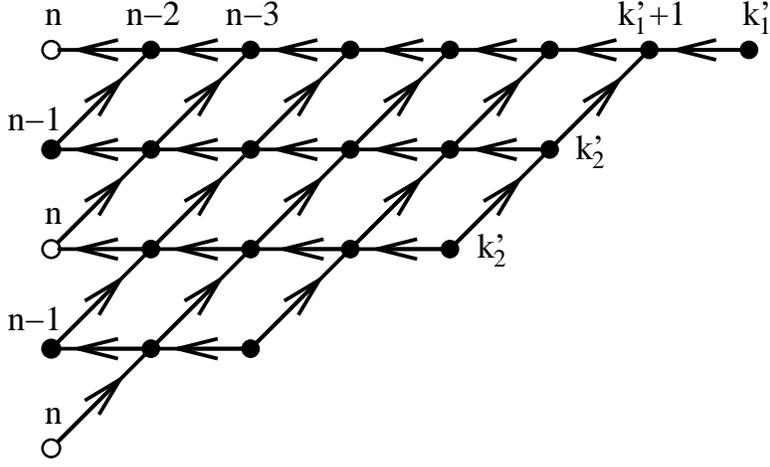,width=4in} \caption{The irreducible
components of $\cup_v \mathcal{L}(\mathbf{v}, \mathbf{w}^n)$ are
enumerated by Young diagrams. \label{fig:youngdiag_dn}}
\end{center}
\end{figure}

The Young diagrams enumerating the irreducible components of
$\cup_\mathbf{v} \mathcal{L}(\mathbf{v}, \mathbf{w}^n)$ can be
visualized as in Figure~\ref{fig:youngdiag_dn}.  For the case of
$\cup_\mathbf{v} \mathcal{L}(\mathbf{v}, \mathbf{w}^{n-1})$,
simply interchange the labellings $n-1$ and $n$. Note that the
vertices in our diagram correspond to the boxes in the classical
Young diagram, and our arrows intersect the classical diagram
edges (cf. \cite{FS03}).

It is relatively easy to compute the geometric action of the
generators $E_k$ and $F_k$ of $\mathfrak{g}$.  We first note that
for every $\mathbf{v}$, $\mathcal{L}(\mathbf{v},\mathbf{w}^n)$ and
$\mathcal{L}(\mathbf{v},\mathbf{w}^{n-1})$ are either empty or a
point.  This can be seen directly or from the dimension formula
for quiver varieties \cite[Cor. 3.12]{N98}. It follows that each
$X_Y^+$ (resp. $X_Y^-$) is equal to
$\mathcal{L}(\mathbf{v},\mathbf{w}^n)$ (resp.
$\mathcal{L}(\mathbf{v},\mathbf{w}^{n-1})$) for some unique
$\mathbf{v}$ which we will denote $\mathbf{v}_Y^+$ (resp.
$\mathbf{v}_Y^-$).

\begin{lem}
The function $g_{X_Y^\pm}$ corresponding to the irreducible component
$X_Y^\pm$ where $Y \in \mathcal{Y}$ is simply $\id_{X_Y^\pm}$, the function on
$X_Y^\pm$ with constant value one.
\end{lem}

\begin{proof}
This is obvious since $X_Y$ is a point.
\end{proof}

\begin{prop}
\label{prop:g_action}
One has $F_k \id_{X_Y^\pm} = \id_{X_{Y'}^\pm}$ where
$\mathbf{v}_{Y'}^\pm = \mathbf{v}_Y^\pm +
\mathbf{e}^k$ if such a $Y'$ exists and $F_k \id_{X_Y^\pm} = 0$
otherwise.  Also,
$E_k \id_{X_Y^\pm} = \id_{X_{Y''}^\pm}$
where $\mathbf{v}_{Y''}^\pm =
\mathbf{v}_Y^\pm - \mathbf{e}^k$ if such a $Y''$ exists and $E_k \id_{X_Y^\pm}
= 0$ otherwise.
\end{prop}

\begin{proof}
We prove the $\mathcal{L}(\mathbf{v},\mathbf{w}^n)$ case.  The case of
$\mathcal{L}(\mathbf{v},\mathbf{w}^{n-1})$ is analogous.
It is clear from the definitions that $F_k \id_{X_Y^+} = c_1
\id_{X^+_{Y'}}$ and $E_k \id_{X_Y^+} = c_2 \id_{X^+_{Y''}}$ for some
constants $c_1$ and $c_2$ if $Y'$ and $Y''$ exist as described above
and that these actions are zero otherwise.  We simply have to compute
the constants $c_1$ and $c_2$. Now,
\begin{align*}
F_k \id_{X_Y^+} (x) &= (\pi_2)_! \pi_1^* \id_{X_Y^+}(x) \\
&= \chi (\{S\, |\, S \text{ is $x$-stable},\, x|_S \in X_Y^+\}) \\
&= \chi (\text{pt}) \\
&= 1
\end{align*}
if $x \in X_{Y'}$ where $\mathbf{v}^+_{Y'} = \mathbf{v}^+_Y +
\mathbf{e}^k$ and zero otherwise.  The fact that the above set is
simply a point follows from the fact that $S_k$ must be the sum of
the images of $x_h$ such that $\inc(h)=k$.  That is, it must be
the span of all the vectors corresponding to the vertices in the
column associated to $S_k$ except the bottommost vertex.  Thus
$c_1=1$ as desired.

A similar argument shows that $c_2=1$.  For a Young diagram $Y$
that contains a removable vertex $k$ (that is, a $Y''$ exists as
described above), there is only one way to extend the unique
representation (up to isomorphism) corresponding to $Y''$ (recall
that $\mathcal{L}(\mathbf{v}_{Y''},\mathbf{w})$ is a point) to a
representation corresponding to $Y$ -- it must be the unique such
representation (up to isomorphism).
\end{proof}

We now compute the weights of the functions corresponding to the
various irreducible components of the quiver variety.  Let
$V=\C^{2n}$ with basis $\{a_1,\dots,a_n,b_1,\dots,b_n\}$ and let
\[
Q: V \times V \to \C
\]
be the nondegenerate, symmetric bilinear form on $V$ given by
\begin{gather*}
Q(a_i,b_i) = Q(b_i,a_i) = 1 \\
Q(a_i,a_j) = Q(b_i,b_j) = 0 \\
Q(a_i,b_j) = 0 \text{ if } i \ne j.
\end{gather*}
Then the Lie algebra $\g = \mathfrak{so}_{2n} \C$ of type $D_n$
consists of the endomorphisms $A : V \to V$ satisfying
\[
Q(Av,w)=Q(v,Aw)=0 \ \forall \ v,w \in V.
\]
In the above basis, the Cartan subalgebra $\mathfrak{h}$ is spanned
by the matrices
\[
D_i = e_{i,i} - e_{n+i,n+i}
\]
where $e_{i,j}$ is the matrix with a one in entry $(i,j)$ and zeroes
everywhere else.  Thus the dual space $\mathfrak{h}^*$ is spanned by
the functions $\{\varepsilon_i\}_{i=1}^n$ given by
\[
\varepsilon_i(D_j) = \delta_{ij}.
\]
The simples roots are given in this basis by
\begin{align*}
\alpha_i &= \varepsilon_i - \varepsilon_{i+1},\ 1\le i \le n-1\\
\alpha_n &= \varepsilon_{n-1} + \varepsilon_n
\end{align*}
and
\begin{align*}
\Lambda_{n-1} &= (\varepsilon_1 + \dots + \varepsilon_{n-1} -
\varepsilon_n)/2 \\
\Lambda_n &= (\varepsilon_1 + \dots + \varepsilon_{n-1} + \varepsilon_n)/2 .
\end{align*}

Let $\lceil a \rceil$ denote the least integer greater than or equal
to $a$ and let $\lfloor a \rfloor$ denote the greatest integer less
than or equal to $a$.

\begin{prop}
\label{prop:weights} For a Young diagram (or partition) $Y =
(\lambda_1, \dots, \lambda_s) \in \mathcal{Y}$, let $\mu = (\mu_1,
\dots, \mu_t)$ be the conjugate partition. Then the weight of
$\id_{X_Y^+}$ is
\begin{gather*}
\Lambda_n - \left\lceil \frac{\mu_1}{2} \right\rceil \alpha_n -
  \left\lfloor
  \frac{\mu_1}{2} \right\rfloor \alpha_{n-1} - \sum_{i=2}^t \mu_i
  \alpha_{n-i} \\
= \frac{1}{2} \left( \sum_{i=1}^{n-1} \varepsilon_i -
\sum_{i=1}^{s}
  \varepsilon_{n-\lambda_i} \right) + \frac{1}{2}
\begin{cases}
\varepsilon_n & \text{if $l(\lambda)$ is even} \\
-\varepsilon_n & \text{if $l(\lambda)$ is odd}
\end{cases}
\end{gather*}
and the weight of $\id_{X_Y^-}$ is
\begin{gather*}
\Lambda_{n-1} - \left\lfloor \frac{\mu_1}{2} \right\rfloor \alpha_n -
  \left\lceil
  \frac{\mu_1}{2} \right\rceil \alpha_{n-1} - \sum_{i=2}^t \mu_i
  \alpha_{n-i} \\
= \frac{1}{2} \left( \sum_{i=1}^{n-1} \varepsilon_i -
\sum_{i=1}^{s}
  \varepsilon_{n-\lambda_i} \right) + \frac{1}{2}
\begin{cases}
\varepsilon_n & \text{if $l(\lambda)$ is odd} \\
-\varepsilon_n & \text{if $l(\lambda)$ is even}
\end{cases}.
\end{gather*}
\end{prop}
\begin{proof}
The first expression follows easily from counting the vertices in the
Young diagram and the second from switching to the bases given by the
$\varepsilon_i$.  For instance, we see that the space $V(k,n-1)$
contributes a weight
\[
-\sum_{i=k}^{n-1} \alpha_i = -\sum_{i=k}^{n-1} (\varepsilon_i -
\varepsilon_{i+1}) = \varepsilon_n - \varepsilon_k
\]
and the space $V(k,n)$, $k\le n-2$, contributes a weight
\[
-\alpha_n -\sum_{i=k}^{n-2} \alpha_i = -(\varepsilon_{n-1} +
\varepsilon_n) - \sum_{i=k}^{n-2} (\varepsilon_i - \varepsilon_{i+1})
= -\varepsilon_n - \varepsilon_k.
\]
\end{proof}


\section{Geometric Realization of the Clifford Algebra}

The standard construction of the spin representations considered above
is through the use of the associated Clifford algebra.  In this
section, we will given a geometric realization of this
algebra.

\subsection{Preliminaries}

We first review the neccessary details concering the Clifford
algebra.  Proofs can be found (for example) in \cite{FH}.

Let the vector space $V$ and the bilinear form $Q$ be as in the
previous section.  Let $C = C(Q)$ be the Clifford algebra associated
to $Q$.  That is, it is the associative algebra generated by $V$ with
relations $v \cdot v = \frac{1}{2}Q(v,v)\cdot 1$ or equivalently
\[
\{v,w\} = v \cdot w + w \cdot v = Q(v,w) \ \forall \ v,w \in V.
\]
Since the relations are linear combinations of elements of even
degree, $C$ inherits a $\Z / 2\Z$-grading $C = C^{\text{even}} \oplus
C^{\text{odd}}$.  Obviously, $C^{\text{even}}$ is a subalgebra of
$C$.

Let $W = \Span \{a_1,\dots,a_n\}$, $W' = \Span\{b_1,\dots,b_n\}$
and $\bigwedge W' = \bigwedge^0 W' \oplus \dots \oplus \bigwedge^n W'$.
Recall that $b_I = b_{i_1} \wedge b_{i_2} \wedge \dots \wedge b_{i_k}$
for $I = \{i_1 < i_2 < \dots < i_k\}$, and with $b_\emptyset=1$ form a
basis for $\bigwedge W'$.

For $w' \in W'$, let $L_{w'} \in \End(\bigwedge W')$ be left
multiplication by $w'$:
\[
L_{w'}(\xi) = w' \wedge \xi,\ \xi \in {\textstyle \bigwedge} W'.
\]
For $\vartheta \in (W')^*$, let $D_\vartheta \in \End(\bigwedge
W')$ be the derivation of $\bigwedge W'$ given by
\[
D_\vartheta (w_1' \wedge \dots \wedge w_r') = \sum_{i=1}^r (-1)^{i-1}
\vartheta(w_i')(w_1' \wedge \dots \wedge {\hat w_i'} \wedge \dots w_r')
\]
where $\hat w_i'$ means the factor $w_i'$ is missing.  Define a
map $l : C \to \End(\bigwedge W')$ by
\[
l(w') = L_{w'}, \quad l(w) = D_\vartheta, \quad w \in W, w' \in W'
\]
where $\vartheta \in (W')^*$ is given by $\vartheta(w') = Q(w,w')$
for all $w' \in W'$.

\begin{lem}
The map $l$ is an isomorphism of algebras $C \cong \End(\bigwedge W')$.
\end{lem}

It follows that
\[
C^{\text{even}} \cong \End({\textstyle \bigwedge}^{\text{even}} W') \oplus
\End({\textstyle \bigwedge}^{\text{odd}} W').
\]
From now on, we suppress the isomophism $l$ and identify $C$ with
$\End(\bigwedge W')$. We can view $C$ as a Lie algebra with bracket
given by the commutator.

\begin{prop}
The map $\g \to C^{\text{even}}$ given on generators by
\begin{align*}
E_k &\mapsto b_{k+1} a_k = L_{b_{k+1}} D_{b_k^*}, \ 1 \le k \le n-1 \\
F_k &\mapsto b_k a_{k+1} = L_{b_k} D_{b_{k+1}^*}, \ 1 \le k \le n-1 \\
E_n &\mapsto a_n a_{n-1} = D_{b_n^*} D_{b_{n-1}^*} \\
F_n &\mapsto b_{n-1} b_n = L_{b_{n-1}} L_{b_n}
\end{align*}
is an injective morphism of Lie algebras.  This gives
$\bigwedge W'$ the
structure of a \g-module and
\[
{\textstyle \bigwedge} W' \cong
L_{\Lambda_{n-1}} \oplus L_{\Lambda_n}
\]
as \g-modules where $L_\lambda$ is the irreducible representation of
highest weight $\lambda$.  In particular
\begin{align*}
{\textstyle \bigwedge}^{\text{even}} W' &\cong
\begin{cases}
L_{\Lambda_n} & \text{ if $n$ is even} \\
L_{\Lambda_{n-1}} & \text{ if $n$ is odd}
\end{cases},\\
{\textstyle \bigwedge}^{\text{odd}} W' &\cong \begin{cases}
L_{\Lambda_n} & \text{ if $n$ is odd} \\
L_{\Lambda_{n-1}} & \text{ if $n$ is even}
\end{cases}.
\end{align*}
The weight of the natural basis vector $b_I$ is
\[
\frac{1}{2}\left( \sum_{i \not \in I} \varepsilon_i - \sum_{j \in I}
  \varepsilon_j \right).
\]
\end{prop}


\subsection{Geometric Realization}

We now present a geometric construction of the Clifford algebra.
Recall the maps
\[
\mathcal{L}(\mathbf{v'}, \mathbf{w}) \stackrel{\pi_1}{\leftarrow}
\mathcal{F}(\mathbf{v}, \mathbf{w}; \mathbf{v} - \mathbf{v'})
\stackrel{\pi_2}{\rightarrow} \mathcal{L}(\mathbf{v}, \mathbf{w})
\]
of Diagram~\ref{eq:diag_action_mod}.  For $1 \le k \le n-1$, define
elements
$\mathbf{y}^k, \mathbf{z}^k \in \Z_{\ge 0}^n$ by
\begin{align*}
\mathbf{y}^k_i &= \begin{cases}
1 & \text{if $k \le i \le n-1$} \\
0 & \text{otherwise}
\end{cases} \\
\mathbf{z}^k_i &= \begin{cases}
1 & \text{if $k \le i \le n-2$ or $k=n$} \\
0 & \text{otherwise}
\end{cases}.
\end{align*}
Then consider the maps
\begin{equation}
\label{eq:cliff_action}
\bigcup_v \mathcal{L}(\mathbf{v}, \mathbf{w}) \stackrel{\pi_1}{\leftarrow}
\bigcup_v \mathcal{F}(\mathbf{v}, \mathbf{w}; \mathbf{y}^k) \cup
\bigcup_v \mathcal{F}(\mathbf{v}, \mathbf{w}; \mathbf{z}^k)
\stackrel{\pi_2}{\rightarrow} \bigcup_v \mathcal{L}(\mathbf{v}, \mathbf{w})
\end{equation}
where $\pi_1$ and $\pi_2$ are as above. Now, for $\mathbf{v} \in
\Z_{\ge 0}^n$, define $\mathbf{v}^*$ by
\begin{gather*}
\mathbf{v}^*_i = \mathbf{v}_i,\ 1 \le i \le n-2 \\
\mathbf{v}^*_{n-1} = \mathbf{v}_n, \quad \mathbf{v}^*_n = \mathbf{v}^*_{n-1}.
\end{gather*}
Note that $\mathbf{v} \mapsto \mathbf{v}^*$ is an involution.  Let
\[
\kappa : \mathcal{L}(\mathbf{v},\mathbf{w}) \to \mathcal{L}(\mathbf{v}^*,
\mathbf{w}^*)
\]
be the map that ``switches'' the vertices $n-1$ and $n$.
Specifically, for $[(x,j)] \in \mathcal{L}(\mathbf{v},\mathbf{w})$
($[(x,j)]$ denotes the $G_\mathbf{V}$-orbit through the point $(x,j)$),
$\kappa([(x,j)]) = [(x',j')]$ where
\begin{align*}
x'_h &= x_h \text{ for } h \ne h_{n-1,n-2}, h_{n,n-2}, h_{n-2,n-1},
h_{n-2,n} \\
x'_{h_{n-1,n-2}} &= x_{h_{n,n-2}},\ x'_{h_{n-2,n-1}} = x_{h_{n-2,n}} \\
x'_{h_{n,n-2}} &= x_{h_{n-1,n-2}},\ x'_{h_{n-2,n}} = x_{h_{n-2,n-1}} \\
j'_i &= j_i \text{ for } i \ne n-1,n \\
j'_{n-1} &= j_n \\
j'_n &= j_{n-1}.
\end{align*}
Note that $\kappa^2$ is the identity map.

Recall that the set
of $\mathbf{v}$ such that $\mathcal{L}(\mathbf{v},\mathbf{w})$ is
non-empty is in 1-1 correspondence with the set $\mathcal{Y}$.
Denote the Young diagram corresponding to such a $\mathbf{v}$ by
$Y_\mathbf{v}$.  Recall that each row in the Young diagram
corresponds to a string of vertices of the Dynkin diagram of \g.  We
say that the \emph{endpoint} of the row is
the lowest index of the vertices appearing in that string with one
exception: we say that the endpoint of the string consisting of the
single vertex $n$ is $n-1$.  For a
Young diagram $Y$, let $l_i(Y)$ be the number
of rows of $Y$ with endpoint strictly less than $i$ (or, equivalently,
the number rows of length greater that $n-i$) and
$l(Y)$ the number of rows of $Y$.  Note that
\[
l(Y_\mathbf{v}) = \mathbf{v}_{n-1} + \mathbf{v}_n.
\]
Now, for $\mathbf{w} = \mathbf{w}^{n-1} \text{ or } \mathbf{w}^n$
and for $1 \le k \le n-1$ define
\begin{align*}
&a_k : M(\mathcal{L}(\mathbf{v},\mathbf{w})) \to
M(\mathcal{L}(\mathbf{v}-\mathbf{y}^k,\mathbf{w^*})) \oplus
M(\mathcal{L}(\mathbf{v-\mathbf{z}^k},\mathbf{w^*})) \\
&b_k : M(\mathcal{L}(\mathbf{v},\mathbf{w})) \to
M(\mathcal{L}(\mathbf{v}+\mathbf{y}^k,\mathbf{w^*})) \oplus
M(\mathcal{L}(\mathbf{v}+\mathbf{z}^k,\mathbf{w^*}))
\end{align*}
by
\begin{align*}
a_k(f) &= (-1)^{l_k(Y_\mathbf{v})} \kappa^* (\pi_1)_! \pi_2^*f \\
b_k(f) &= (-1)^{l_k(Y_\mathbf{v})} \kappa^* (\pi_2)_! \pi_1^*f
\end{align*}
where we have used Diagram~\ref{eq:cliff_action}.
Define
\begin{align*}
&a_n : M(\mathcal{L}(\mathbf{v},\mathbf{w})) \to
M(\mathcal{L}(\mathbf{v},\mathbf{w^*})) \\
&b_n : M(\mathcal{L}(\mathbf{v},\mathbf{w})) \to
M(\mathcal{L}(\mathbf{v},\mathbf{w^*}))
\end{align*}
by
\begin{align*}
a_n(f) &= \begin{cases}
(-1)^{l(Y_\mathbf{v})} \kappa^* f & \text{if $\dim \mathbf{w}_n +
  l(Y_\mathbf{v})$ is even} \\
0 & \text{if $\dim \mathbf{w}_n + l(Y_\mathbf{v})$ is odd}.
\end{cases}\\
b_n(f) &= \begin{cases}
(-1)^{l(Y_\mathbf{v})} \kappa^* f & \text{if $\dim \mathbf{w}_n +
  l(Y_\mathbf{v})$ is odd} \\
0 & \text{if $\dim \mathbf{w}_n + l(Y_\mathbf{v})$ is even}.
\end{cases}
\end{align*}

For $\mathbf{v} = \mathbf{v}_Y^\pm$ for some $Y \in \mathcal{Y}$, note
that there exists a $Y^{+k}$ such that $\mathbf{v}^\pm_{Y^{+k}} =
\mathbf{v}^\pm_Y +
\mathbf{y}^k$ or $\mathbf{v}^\pm_Y + \mathbf{z}^k$ if and only
if $Y$ does not contain a row with endpoint $k$.  If this is the case,
then $Y^{+k}$ is obtained from $Y$ by adding such a row and if
$\mathbf{v}^\pm_{Y^{+k}} = \mathbf{v}^\pm_Y + \mathbf{y}^k$, then
$\mathcal{L}(\mathbf{v}^\pm_Y + \mathbf{z}^k)$ is empty and vice
versa.  Similarly, there exists a $Y^{-k}$ such that
$\mathbf{v}^\pm_{Y^{-k}} = \mathbf{v}^\pm_Y -
\mathbf{y}^k$ or $\mathbf{v}^\pm_Y - \mathbf{z}^k$ if and only
if $Y$ contains a row with endpoint $k$.  If this is the case,
then $Y^{-k}$ is obtained from $Y$ by removing such a row and if
$\mathbf{v}^\pm_{Y^{-k}} = \mathbf{v}^\pm_Y - \mathbf{y}^k$, then
$\mathcal{L}(\mathbf{v}^\pm_Y - \mathbf{z}^k)$ is empty and vice
versa.

\begin{prop}
\label{prop:kadd_action}
Using Diagram~\ref{eq:cliff_action},
\begin{align*}
(\pi_1)_!\pi_2^* \id_{X^\pm_Y} &= \begin{cases}
\id_{X^\pm_{Y^{-k}}} & \text{if $Y$ contains a row with endpoint $k$} \\
0 & \text{otherwise}
\end{cases}, \\
(\pi_2)_!\pi_1^* \id_{X^\pm_Y} &= \begin{cases}
\id_{X^\pm_{Y^{+k}}} & \text{if $Y$ does not contain a row with endpoint $k$} \\
0 & \text{otherwise}
\end{cases}, \\
\kappa^* \id_{X^\pm_Y} &= \id_{X^\mp_Y}.
\end{align*}
\end{prop}

\begin{proof}
It is clear from the definitions that $(\pi_1)_!\pi_2^* \id_{X^\pm_Y} =
c_{-k} \id_{X^\pm_{Y^{-k}}}$ for some constant $c_{-k}$ if $Y$ contains a
row with endpoint $k$ and is zero otherwise.  Similarly, it is clear that
$(\pi_2)_!\pi_1^* \id_{X^\pm_Y} =
c_{+k} \id_{X^\pm_{Y^{+k}}}$ for some constant $c_{+k}$ if $Y$ does
not contain a
row with endpoint $k$ and is zero otherwise.  It remains to compute the
constants $c_{-k}$ and $c_{+k}$.  Assume $Y$ does not contain a row with endpoint
$k$. Then,
\begin{align*}
(\pi_2)_!\pi_1^* \id_{X^\pm_Y} (x)
&= \chi(\{S\ |\ S \text{ is $x$-stable},\ x|_S \in X^\pm_Y\}) \\
&= \chi(\text{pt}) \\
&= 1
\end{align*}
if $x \in X^\pm_{Y^{+k}}$ and is zero otherwise.  The fact that the
above set is a point follows from the fact that $x$-stability implies
that $S$ must be the span of the vectors corresponding to the vertices of
the Young diagram obtained from $Y^{+k}$ by removing the lowest vertex
of each of the first $n-k$ columns.  This is because the vector
corresponding to a given vertex lies in the smallest subspace of
$\mathbf{V}$ containing any vector below it, in the same column.  Thus
$c_{+k}=1$.

A similar argument shows that $c_{-k}=1$.  For a Young diagram $Y$ that
contains a row with endpoint $k$, there is only one way to
extend the unique representation (up to isomorphism) corresponding to
$Y^{-k}$ (recall that
$\mathcal{L}(\mathbf{v}_Y,\mathbf{w})$ is a point) to a representation
corresponding to $Y$ -- it must be the unique such representation (up
to isomorphism).
\end{proof}

\begin{theo}
The operators $a_k$, $b_k$ on $L(\mathbf{w}^{n-1}) \oplus
L(\mathbf{w}^n)$ provide the structure of a representation of the
Clifford algebra $C$.  Moreover, as representations of $C$,
\begin{equation}
\label{eq:cliff_isom}
L(\mathbf{w}^{n-1}) \oplus L(\mathbf{w}^n) \cong {\textstyle \bigwedge} W.
\end{equation}
\end{theo}

\begin{proof}
For $Y = (n-1 \ge l_1 > \dots > l_s \ge 1) \in \mathcal{Y}$, let
\begin{align*}
I_Y^+ &= \begin{cases}
\{l_s,\dots,l_1\} & \text{if $s$ is even} \\
\{l_s,\dots,l_1,n\} & \text{if $s$ is odd}
\end{cases} \\
I_Y^- &= \begin{cases}
\{l_s,\dots,l_1\} & \text{if $s$ is odd} \\
\{l_s,\dots,l_1,n\} & \text{if $s$ is even}
\end{cases}.
\end{align*}
Both sides of \eqref{eq:cliff_isom} have dimension $2^n$. We
identify the two via the map $\id_{X^\pm_Y} \mapsto b_{I_Y^\pm}$.
For a set $I$ of indices between $1$ and $n$ that does not contain
the index $k$, let $I^{+k}$ be the set obtained by adding $k$. For
a set $I$ of indices between $1$ and $n$ that contains the index
$k$, let $I^{-k}$ be the set obtained by removing $k$.  We view
such sets of indices as a partition or Young diagram as follows:
for $I = \{i_1,\dots,i_s\}$, let the corresponding Young diagram
have rows of length $n-i_j$ (with endpoint $i_j$) for those $1 \le
j \le s$ such that $i_j < n$. Thus $l_i(I)$ is defined -- it is
the number of indices of $I$ strictly less than $i$. Now, for $1
\le k \le n$,
\begin{align*}
a_k (b_{I_Y^\pm}) &= D_{b_k^*} (b_{I_Y^\pm})
= \begin{cases}
(-1)^{l_k(I_Y)} b_{(I_Y^\pm)^{-k}} & \text{if $I_Y^\pm$ contains $k$} \\
0 & \text{if $I_Y^\pm$ does not contain $k$}
\end{cases}, \\
b_k (b_{I_Y^\pm}) &= L_{b_k} (b_{I_Y^\pm})
= \begin{cases}
(-1)^{l_k(I_Y)} b_{(I_Y^\pm)^{+k}} & \text{if $I_Y^\pm$ does not
  contain $k$} \\
0 & \text{if $I_Y^\pm$ contains $k$}
\end{cases}.
\end{align*}
Note that $(I_Y^\pm)^{\pm k} = I_{Y^{\pm k}}^\mp$ for $1 \le k \le
n-1$, $(I_Y^\pm)^{+ n} = I_Y^\mp$ if $I_Y^\pm$ does not contain
$n$ and $(I_Y^\pm)^{- n} = I_Y^\mp$ if $I_Y^\pm$ contains $n$.

Now, recalling that a plus sign indicates the case $\mathbf{w} =
\mathbf{w}^n$ and a minus sign indicates the case $\mathbf{w} =
\mathbf{w}^{n-1}$, we see that $I_Y^\pm$ contains an $n$ if and only
if $\mathbf{w}_n + l(Y) = \mathbf{w}_n + l_n(Y)$ is even.  Thus the
result follows from the definition of the geometric action of the
$a_k$ and $b_k$ and Proposition~\ref{prop:kadd_action}.
\end{proof}

We now show that the geometric action of the Clifford algebra we have
constructed is a natural extension of the geometric action of \g.

\begin{prop}
As operators on $L(\mathbf{w}^{n-1}) \oplus L(\mathbf{w}^n)$,
\begin{align*}
E_k &= b_{k+1} a_k, \quad 1 \le k \le n-1 \\
F_k &= b_k a_{k+1}, \quad 1 \le k \le n-1 \\
E_n &= a_n a_{n-1} \\
F_n &= b_{n-1} b_n.
\end{align*}
\end{prop}
\begin{proof}
It suffices to show that the relations hold on the $\id_{X^\pm_Y}$
since these functions span $L(\mathbf{w}^{n-1}) \oplus
L(\mathbf{w}^n)$.
Note that we can remove a vertex with degree $k$ from $Y$ (and be left
with an element of $\mathcal{Y}$) if and only
if it has a row with endpoint $k$ but no row with endpoint $k+1$.  Let
$Y$ be such a diagram.
For $1 \le k \le n-1$,
\begin{align*}
b_{k+1} a_k \id_{X^\pm_Y}
&= (-1)^{l_k(Y)} b_{k+1} \kappa^* (\pi_1)_! \pi_2^* \id_{X^\pm_Y} \\
&= (-1)^{l_k(Y)} b_{k+1} \kappa^* \id_{X^\pm_{Y^{-k}}} \\
&= (-1)^{l_k(Y)} b_{k+1} \id_{X^\mp_{Y^{-k}}} \\
&= (-1)^{l_k(Y)} (-1)^{l_{k+1}(Y^{-k})} \kappa^* (\pi_2)^! \pi_1^*
\id_{X^\mp_{Y^{-k}}} \\
&= (-1)^{l_k(Y)} (-1)^{l_{k+1}(Y^{-k})} \kappa^*
\id_{X^\mp_{(Y^{-k})^{+(k+1)}}} \\
&= \id_{X^\pm_{(Y^{-k})^{+(k+1)}}}
\end{align*}
But $(Y^{-k})^{+(k+1)}$ is precisely the Young diagram obtained from
$Y$ by removing a vertex $k$ and so $b_{k+1} a_k \id_{X^\pm_Y} = E_k
\id_{X^\pm_Y}$.  The proof of the other three equations is analogous.
\end{proof}


\section{Extensions}
The results of this paper can be easily extended to type
$D_\infty$. One obtains an enumeration of the irreducible
components of Nakajima's quiver variety by Young diagrams of
stricly decreasing row lengths but with no condition on maximum
row length.  The geometric Lie algebra and Clifford algebra
actions can be easily computed and are analogous to those obtained
for type $D_n$.

If one is interested in the $q$-deformed analogues of the Clifford
algebra and spin representations (cf. \cite{DF94}), this can also
be handled by the methods of the current paper.  One simply has to
consider the quiver varieties over a finite field instead of over
$\C$.

To treat arbitrary finite dimensional irreducible representations
instead of just the spin representations considered here, somewhat
different techniques are required.  This situation will be examined
in a upcoming paper.


\small{
}
\vspace{4mm}
\noindent
Alistair Savage,\\
The Fields Institute for Research in Mathematical Sciences and
University of Toronto, Toronto, Ontario, Canada,\\
email:\;\texttt{alistair.savage@aya.yale.edu}

\end{document}